\documentclass[notitlepage,leqno,10pt]{article}
\usepackage{latexsym}
\usepackage{amsmath}
\usepackage{amsfonts}
\usepackage{amssymb}

\renewcommand{\a }{\alpha }

\renewcommand{\d}{\delta }
\newcommand{\D }{\Delta }

\newcommand{\e }{\varepsilon }
\newcommand{\g }{\gamma}

\renewcommand{\l }{\lambda }

\newcommand{\n }{\nabla }
\newcommand{\var }{\varphi }

\newcommand{\s }{\sigma }
\newcommand{\Sig }{\Sigma}

\newcommand{\ov}{\overline}
\newcommand{\intbar}{\mathop{\int\makebox(-13.5,0){\rule[4pt]{.7em}{0.3pt}}%
\kern-6pt}\nolimits}

\newcommand{\be}{\begin{equation}}
\newcommand{\ee}{\end{equation}}
\newcommand{\bes}{\begin{equation*}}
\newcommand{\ees}{\end{equation*}}
\newcommand{\ba}{\begin{eqnarray}}
\newcommand{\ea}{\end{eqnarray}}
\newcommand{\bas}{\begin{eqnarray*}}
\newcommand{\eas}{\end{eqnarray*}}
\newenvironment{pf}{\noindent{\sc Proof}.\enspace}{\rule{2mm}{2mm}\medskip}
\newenvironment{pfn}{\noindent{\sc Proof}}{\rule{2mm}{2mm}\medskip}

\newcommand{\R}{\mathbb{R}}

\newcommand{\Z}{\mathbb{Z}}

\newcommand{\N}{\mathbb{N}}

\author{ Cheikh Birahim NDIAYE}

\date{}

\title{ \bf Existence results for mean field equation with turbulence}

\begin{document}

\newtheorem{lem}{Lemma}[section]
\newtheorem{pro}[lem]{Proposition}
\newtheorem{thm}[lem]{Theorem}
\newtheorem{rem}[lem]{Remark}
\newtheorem{cor}[lem]{Corollary}
\newtheorem{df}[lem]{Definition}

\maketitle

\begin{center}

{\small SISSA, via Beirut 2-4, 34014 Trieste, Italy.}

\end{center}

\

\
\begin{center} 
{\bf Abstract}
\end{center}
In this paper we consider the following form of the so-called Mean field equation arising from the statistical mechanics description of two dimensional turbulence
\begin{equation}\label{eq:study}
- \D_g u = \rho_1
 \left( \frac{ e^{u}}{\int_\Sig e^{u} dV_g}-1\right)-\rho_2
 \left( \frac{ e^{-u}}{\int_\Sig e^{-u} dV_g}
- 1 \right)
\end{equation} 
on a given closed orientable Riemannian surface \;($\Sigma,\;g$) with volume \;$1$, where $\rho_1, \rho_2$ are real parameters. Exploiting the variational structure of the problem and running a min-max scheme introduced by Djadli and Malchiodi, we prove that if  \;$k$\;is a positive integer, \;$\rho_1$\;and \;$\rho_2$\;two real numbers such that \;$\rho_1\in (8k\pi, 8(k+1)\pi)$\;and\;$\rho_2<4\pi$\;then \;$\eqref{eq:study}$\;is solvable.

\begin{center}

\bigskip\bigskip


\bigskip

\centerline{\bf AMS subject classification: 35B33, 35J35, 53A30,
53C21}

\end{center}

\footnotetext[1]{E-mail addresses: ndiaye@sissa}

\section{Introduction}
Many problems in physics can be formulated in terms of nonlinear elliptic equations with exponential nolinearities.\\
A typical example is the so called {\em mean field equation} on a given closed Riemannian surface ($\Sigma,g$)  with volume \;$1$.
\begin{equation}\label{eq:meanf}
- \D_g u = \rho
 \left( \frac{ he^{u}}{\int_\Sig he^{u} dV_g}-1\right)\;\;\text{on}\;\Sigma;
\end{equation}
( where \;$\D_g$ is the Laplace-Beltrami, $\rho$\;a real parameter) which arises in the study of limit of point vortices of Euler flows,
spherical Onsager vortex theory and condensates in some
Chern-Simons-Higgs models, see for example the papers \cite{bt},
\cite{bclt}, \cite{clmp1}, \cite{clmp2}, \cite{cl1}, \cite{cl2},
\cite{djlw}, \cite{ki}, \cite{st} and the references therein.\\\\
An other example is what we refer as  {\em mean field equation with turbulence} on a closed orientable Riemannian surface \;($\Sigma,g$)
 \begin{equation}\label{eq:meanft}
- \D_g u = \rho_1
 \left( \frac{ e^{u}}{\int_\Sig e^{u} dV_g}-\frac{1}{|\Sigma|}\right)-\rho_2
 \left( \frac{ e^{-u}}{\int_\Sig e^{-u} dV_g}
- \frac{1}{|\Sigma|} \right);\;\;\;\int_{\Sigma}udV_g=0.
\end{equation}
(where \;$\D_g$ is the Laplace-Beltrami,\;$|\Sigma|$\;the volume of \;$\Sigma$ and $\rho_1$ and \;$\rho_2$\;are two positive real parameters)
arising from the statistical mechanics description of two-dimensional turbulence see Joyce and Montgomery \cite{jm} and Pointin and Lundgren \cite{pl}.\\

The {\em mean field equation}\;$\eqref{eq:meanf}$\; has received much attention in the last two decades. To mention some related non-trivial results, we cite the one of Ding-Jost-Li-Wang which asserts that if the underlying surface \;$\Sig$\;\;has positive genus then the problem has a solution provided \;$\rho\in]8\pi,16\pi[$, see \cite{djlw}. Latter, using degree theory argument, Chen and Lin improve Ding-Jost-Li-Wang result by showing that if still the genus is positive then the problem is solvable for every \;$\rho\neq k8\pi$\;\;where \;$k$\;is an arbitrary positive integer, see \cite{cl}. Recently, Zindine Djadli refines Chen and Lin result by removing the constraint on the genus, see \cite {dj}. In the critical case, namely when \;$\rho=8\pi$, Ding-Jost-Li-Wang have given sufficient conditions for the solvability. \\

From this panorama on the {\em mean field equation}, we see that the answer to the question of existence of solutions is quite satisfactory. However for the {\em mean field equation with turbulence}\;$\eqref{eq:meanft}$, little is known. To the best of our knowledge, the only available result in the literature is the one of Ohtsuka-Suzuki\cite{os} and Ricciardi\cite{rt}. In fact Ohtsuka-Suzuki obtained existence of solutions for \;$\rho_i\in[0,8\pi[$\;via minimization and Ricciardi prove recently existence of Mountain-pass solutions under the assumptions that\;($\Sig,g$) is a closed Riemannian surface\;such that the first non-zero eigenvalue \;$\mu_1(\Sig)$\;of \;$-\D_g$\; verifies\;$8\pi<\mu_1(\Sig)|\Sig|<16\pi$\; and for \;$\rho_1$\;and\;$\rho_2$\;such that\;$\rho_1+\rho_2<\mu_1(\Sig)|\Sig|$\;and \;$\max_{i=1,2}(\rho_i)>8\pi$.\\

In this paper we will consider the following version of the {\em mean field equation with turbulence}
\begin{equation}\label{eq:eqs}
- \D_g u = \rho_1
 \left( \frac{ e^{u}}{\int_\Sig e^{u} dV_g}-1\right)-\rho_2
 \left( \frac{ e^{-u}}{\int_\Sig e^{-u} dV_g}
- 1 \right).
\end{equation}
where \;$\Sigma$\;has volume \;$1$ and the parameters \;$\rho_i$\;are arbitrary real numbers (we recall that the relevant case for physics is when both are non-negative).\\
Prolem\;$\eqref{eq:eqs}$\;is variational. Indeed critical points of the following functional 
\begin{equation}
\text{II}_{\rho}(u)=\frac{1}{2}\int_{\Sigma}|\n u|^2dV_g-\rho_1\log\int_{\Sigma}e^{u-\bar u}dV_g-\rho_2\log\int_{\Sigma}e^{-u+\bar u}dV_g,\;\;\;u\in H^1(\Sigma);
\end{equation}
(where \;$\rho=(\rho_1,\rho_2)$) are weak solutions, hence due to standard elliptic regularity are also classical solutions.
\vspace{6pt}

\indent
Our main goal is to give a more general existence result of the type of  Zindine Djadli for the {\em mean field equation}.

We have indeed the following theorem.
\begin{thm}\label{eq:theo}
Suppose $k$ is a positive integer. Assume that\; $\rho_1\in (k8\pi ,8(k+1)\pi)$\;and  $\rho_2 <4\pi$, then problem\;$\eqref{eq:eqs}$\; is
solvable
\end{thm}
We are going to describe the main ideas in the proof of Theorem\;$\ref{eq:theo}$. From Theorem\;$\ref{th:os1}$\;below, if one of the \;$\rho_i$'s is bigger the \;$8\pi$, then the functional \;$\text{II}_{\rho}$\;is not bounded from below, hence extremals have to be found amongs saddle points. To do so we will use a min-max scheme introduced by Djadli and Malchiodi in their study of the existence of  constant \;$Q$-curvature metrics on four manifolds, see \cite{dm}. By classical arguments in critical point theory, such a scheme yields  existence of {\em Palais-smale} sequences, namely sequences \;$(u_l)$, such that
$$
\text{II}_{\rho}(u_l)\rightarrow c\in \R,\;\;\;\;\text{II}^{'}_{\rho}(u_l)\rightarrow 0.
$$
Since the functional \;$\text{II}_{\rho}(\cdot)$\;is invariant under translation by constant, then we can always assume that the sequence \;$(u_l)$\;satisfies the normalization
$$
\int_{\Sig}e^{u_l}dV_g=1\;\;\;\forall l.
$$
If one proves that such sequences are bounded or that a similar compactness criterion holds, then the existence of solutions to problem\;$\eqref{eq:eqs}$\; follows automatically.\\
To do so, we apply Struwe monotonicity method, see \cite{str}. This consits in performing a min-max argument for
for different values of \;$\rho$\;of the form \;$\rho=(t\rho_1,t\rho_2)$\;\;$t\sim1$, and then to prove that there exists bounded Palais-Smale sequences for\;$\text{II}_{\rho}$\;for \;$\rho=(t_l\rho_1,t_l\rho_2)$\;with \;$t_l\rightarrow1$. This yields existence of solutions to the problems.
\begin{equation*}
- \D_g u = t_l\rho_1
 \left( \frac{ e^{u}}{\int_\Sig e^{u} dV_g}-1\right)-t_l\rho_2
 \left( \frac{ e^{-u}}{\int_\Sig e^{-u} dV_g}
- 1 \right).
\end{equation*}
Hence an application of Proposition\;$\ref{eq:compactness}$, gives the existence of solutions to problem\;$\eqref{eq:eqs}$.\\
From the discussion above, we have that the core of the analysis consist in finding Palais-Smale sequences. This will be done by characterizing the topology of low sublevels of the functional \;$\text{II}_{\rho}$.
From considerations coming from an improvement of the Moser-Trudinger type inequality (Theorem\;$\ref{th:os1}$), it follows that if \;$\text{II}_{\rho}(u)$\;attains large negative values then \;$e^{u}$ has to concentrate near at most \;$k$\;points of\;$\Sig$. This means that, if we normalize \;$u$\;so that \;\;$\int_{\Sig}e^{u}dV_g=1$, then naively \;$e^{u}\simeq\sum_{i=1}^kt_i\d_{x_i},\;\;x_i\in \Sig,\;\; t_i\geq 0,\;\sum_{i=1}^kt_i=1$. Such a family of convex combination of Dirac deltas are called formal barycenters of \;$\Sig$\;of order \;$k$, see Section 2, and will be denoted by \;$\Sig_k$. With a further analysis (see Subsection 3.3 ), it is possible to show that the sublevel\;$\{\text{II}_{\rho}<-L\}$\;for large \;$L$\; has the same homology as \;$\Sig_k$.
Using the non-contractibility of \;$\Sig_k$, we perform a min-max scheme, and get the Palais-Smale sequences.

\vspace{10pt}

\noindent {\bf Acknowledgements: } The author have been supported
by M.U.R.S.T within the PRIN 2006 {\em Variational methods and
nonlinear differential equations}.

\section{Notation and preliminaries}\label{s:pr}

In this section we collect some useful preliminary facts. For $x, y
\in \Sig$ we denote by $d(x,y)$ the metric distance between $x$ and
$y$ on $\Sig$. In the same way, we denote by $d(S_1, S_2)$ the
distance between two sets $S_1, S_2 \subseteq \Sig$, namely
$$
  d(S_1,S_2) = \inf \left\{ d(x,y) \; : \; x \in S_1, y \in S_2
\right\}.
$$
Recalling that we are assuming $Vol_g(\Sig) := \int_\Sig 1 dV_g =
1$, given a function $u \in L^1(\Sig)$, we denote its average (or
integral) as
$$
  \ov{u} = \int_\Sig u dV_g.
$$
Below, by $C$ we denote large constants which are allowed to vary
among different formulas or even within lines. When we want to
stress the dependence of the constants on some parameter (or
parameters), we add subscripts to $C$, as $C_\d$, etc.. Also
constants with subscripts are allowed to vary.

We now recall some Moser-Trudinger type inequalities and compactness
results. The Euler-Langrange functional under study is the following
\begin{equation}
\text{II}_{\rho}(u)=\frac{1}{2}\int_{\Sigma}|\n u|^2dV_g-\rho_1\log\int_{\Sigma}e^{u-\bar u}dV_g-\rho_2\log\int_{\Sigma}e^{-u+\bar u}dV_g,\;\;\;u\in H^1(\Sigma);
\end{equation}
which for large values of $\rho_1$ and $\rho_2$ will be in general
unbounded from below. In fact, there is a precise criterion for
$\text{II}_\rho$ to be bounded from below, which has been proved by Ohtsuka and Suzuki.

\begin{thm}\label{th:os1} (\cite{os})
For $\rho = (\rho_1, \rho_2)$ the functional $\text{II}_\rho 
$ is bounded from below if and only if both $\rho_1$
and $\rho_2$ satisfy the inequality $\rho_i \leq 8 \pi$.
\end{thm}
Now we recall the following classical Moser-Trudinger inequality.
\begin{lem}
There exists a constant \;$C_1$\;depending only on \;$(\Sig,g)$, such that
$$
\int_{\Sig}e^{4\pi(u-\ov{ u})}dV_g\leq C_1\;\;\;\;\forall u\in H^1(\Sig)\;\;\text{such that}\;\;\int_{\Sig}|\n u|^2dV_g\leq 1.
$$
As a consequence we have 
\begin{equation}\label{eq:mt}
    \log \int_\Sig e^{(u - \ov{u})} dV_g \leq C + \frac{1}{16 \pi}
\int_{\Sig} |\n u|^2 dV_g,\;\;\;\forall u\in H^1(\Sig).
\end{equation}
\end{lem}
Next we give a compactness result due to Ohtsuka and Suzuki.
\begin{thm}\label{th:os2}(\cite{os})
Let \;$\rho_{1,n}$\;and \;$\rho_{2,n}$\;be sequences of non-negative real numbers satisfying
$$
\rho_{i,n}\rightarrow \rho_i\;\;\text{as}\;\;n\rightarrow +\infty;
$$
and \;$u_{n}$\;be a sequence of solutions to \;$\eqref{eq:eqs}$\;corresponding to \;$(\rho_{i,n}, \rho_{2,n})$, with \;$\ov{u_{n}}=0$. let also \;$\mu_{i,n}$\;be the following Radon measures 
$$
\mu_{1,n}=\frac{\rho_{1,n}e^{u_{n}}}{\int_{\Sigma}e^{u_n}dV_g}dV_g;
$$
$$
\mu_{2,n}=\frac{\rho_{2,n}e^{-u_{n}}}{\int_{\Sigma}e^{-u_n}dV_g}dV_g.
$$
Moreover let \;$w_{i,n}$\;be as follows
$$
w_{i,n}=\int_{\Sigma}G(x,y)d\mu_{i,n}(y);
$$
where \;$G$\;is the Green function of \;$-\D_g$\;such that \;$\int_{\Sigma}G(\cdot,y)dV_g(y)=0$, and we assume also without loss of generality that
$$
\mu_{i,n}\rightharpoonup \mu_i \;\text{weakly*}
$$
Let \;$S_1$\;and \;$S_2$\;denotes the following sets
$$
S_1=\{x\in \Sigma:\;\exists x_n\in \Sigma\;\;\text{s.t}\;\;u_n(x_n)\rightarrow+\infty\};
$$
and 
$$
S_2=\{x\in \Sigma:\;\exists x_n\in \Sigma\;\;\text{s.t}\;\;u_n(x_n)\rightarrow-\infty\}.
$$
Then the following alternatives hold:\\\\
$(1)$\;( Compactness)\\
We have that \;$S_1\cup S_2=\emptyset$\; and there exists \;$u\in H^1(\Sig)$,\;$\ov{ u}=0$\;and (up to subsequence)\; 
$$
u_n\rightarrow u\;\;\text{in}\;\;H^1(\Sigma)
$$
and \;$u$\;is a solution of \;$\eqref{eq:eqs}$\;for \;$\rho_1$\;and \;$\rho_2$.\\\\
$(2)$\;( one-sided concentration)\\
There exists \;$i\in \{1,2\}$\;such that \;$S_i\neq \emptyset$\;and \;$S_j=\emptyset$\;for \;$j\in \{1,2\}\setminus\{i\}$. Moreover, it holds that
$$
\mu_i=\sum_{x_0\in S_i}8\pi\delta_{x_0}
$$
and
$$
\mu_{i,n}\rightarrow 0\;\text{in}\;\; L^{\infty}(\omega);
$$
for every \;$\omega\subset\subset \Sig\setminus S_i$. On the other hand, there exists \;$w_{j}\in H^1(\Sig)$\;with \;$\ov{ w_j}=0$\; such that up to a subsequence
$$
w_{j,n}\rightarrow w_j\;\;\text{in}\;\;H^1(\Sig)
$$
and \;$w_j$\;is solution to
$$
-\D_gw=\l\left(\frac{Ke^w}{\int_{\Sig}Ke^wdV_g}-1\right),\;\;\;\ov{ v}=0;
$$
with \;$K(x)=e^{-\sum_{x_0\in S_i}8\pi G(x,x_0)}$. \\
$(3)$\;(concentration)\\
For each \;$i=1,2$, we 
have \;$S_i\neq \emptyset$\;and there exists a positive constant \;$m_i(x_0)\geq 4\pi$\;for each \;$x_0\in S_i$. Furthermore, we have a non-negative function \;$r_i\in L^1(\Sig)\cap L^{\infty}_{loc}(M\setminus S_i)$\;such that 
$$
\mu_i=r_i+\sum_{x_0\in S_i}m_i(x_0)\delta_{x_0};
$$
and
$$
\mu_{i,n}\rightarrow r_i\;\;\text{in}\;\;L^p(\omega)
$$
for every \;$p\in [0,+\infty[$\;and every \;$\omega\subset\subset \Sigma \setminus S_i$. Finally the following fact hold:\\
$3-i$)\\
If there exists \;$x_0\in S_i\setminus S_j$\;for \;$i\neq j$, then we have \;$m_i(x_0)=8\pi$ and \;$r_i=0$.\\
$3-ii$)\\
For every \;$x_0\in S_1\cap S_2$, we have
$$
(m_1(x_0)-m_2(x_0))^2=8\pi(m_1(x_0)+m_2(x_0)).
$$
moreover, if \;$S_i\subset S_j$\;and there exists \;$x_0\in S_i$\;satisfying
$$m_i(x_0)-m_j(x_0)>4\pi,$$
then we have \;$r_i=0$.
\end{thm}
Now we recall a Theorem due to Yanyan Li, which will be used  to derive a compactness result adapted to our purposes.
\begin{thm}\label{th:yyli} (\cite{li})
Let $(u_n)_n$ be a sequence of solutions of the equations
$$
  - \D u_n = \l_n \left( \frac{V_n e^{u_n}}{\int_\Sig V_n e^{u_n} dV_g} -
W_n
\right),
$$
where $(V_n)_n$ and $(W_n)_n$ satisfy
$$
  \int_\Sig W_n dV_g = 1; \qquad \|W_n\|_{C^1(\Sig)} \leq C; \qquad |\log V_n| \leq C; \qquad \|\n
V_k\|_{L^\infty(\Sig)} \leq C,
$$
and where $\l_n \to \l_0 > 0$, $\l_0 \neq 8 q \pi$ for $q = 1, 2,
\dots$. Then, under the additional constraint $\int_\Sig u_n dV_g =
0$, $(u_n)_n$ stays uniformly bounded in $L^\infty(\Sig)$.
\end{thm}

Next we give a compactness result which describe all the possibles cases of Theorem\;$\ref{eq:theo}$.
\begin{pro}\label{eq:compactness}
Let \;$K_1$\;be a compact set of \;$\cup_{i=1}^{\infty}(8\pi i, 8\pi(i+1))$\;and \;$K_2$\;be a compact set of \;$(-\infty,4\pi)$. Let \;$\rho_{1,n}$\;be a sequence in \;$K_1$\;and \;$\rho_{2,n}$\;be a sequence in \;$K_2$. Moreover let \;$u_n$\;be a sequence of solutions to \;$\eqref{eq:eqs}$\;correspnding to \;$\rho_{1,n}$\;and\;$\rho_{2,n}$\;with \;$\bar u_n=0$. Then we have \;$u_n$\;is bounded in\;$C^m(\Sig)$\;for every positive integer \;$m$.
\end{pro}
\begin{pf}
We first claim that for every \;$p>1$\;there exists \;$\bar \rho$\;(depending on \;$K_1$,\;$K_2$\;and \;$p$\;)such that

\begin{equation}
\int_{\Sig}e^{-pu_n}dV_g\leq C.
\end{equation}
To prove the Claim we use the Green representation formula for \;$-u_n$, an argument of Brezis and Merle, see \cite{bm} and thhe fact that \;$\rho_{1,n}>0$. Indeed we have that
$$
  -u_{n}(x) \leq C + \int_\Sig G(x,y) \left(2 \rho_{2,n} \frac{e^{-u_{n}}}{\int_\Sig
     e^{-u_{n}} dV_g} \right) dV_g(y),
$$
where $G(x,y)$ is the Green's function of $- \D_g$ on $\Sig$. Next using
Jensen's inequality we find
$$
  e^{-p u_{n}(x)} \leq C \int_\Sig \exp(-2 p \rho_{2,n} G(x,y)) \frac{
  e^{-u_{n}}}{\int_\Sig e^{-u_{n}} dV_g} dV_g(y).
$$
Now using the asymptotics of the Green function ($G(x,y) \simeq \frac{1}{2\pi} \log\left(
\frac{1}{d(x,y)} \right)$) and also the Fubini theorem we get
$$
  \int_\Sig e^{-p u_{n}} dV_g \leq C \sup_{x \in \Sig} \int_\Sig
\frac{1}{d(x,y)^{\frac{p \rho_{2,n}}{\pi}}} dV_g(y).
$$
Thus it is sufficient to take $\ov{\rho} = \frac{\pi}{2 p}$ in order
to obtain the claim.\\
Now suppose \;$\rho_{2,n}\geq \bar \rho$\;. Using Theorem\;$\ref{th:os2}$\;we have that three alternatives can occur. On the other hand since \;$\rho_{i,n}\in K_1$\;and \;$\rho_{2,n}$, then it is trivially seen that  the one-sided concentration and the concentration alternatives can not occur. Hence we have compactness and using standard elliptic regularity theory, we have boundedness in \;$C^{m}(\Sig)$\;for every \;$m$.\\
Now supose\;$\rho_{2,n}\leq \bar \rho$. Then from the Claim, we have \;$e^{-u_n}$\;uniformly bounded in \;$L^p$. Hence \;$v_{n}$\;defined as follows
\begin{equation}
-\D_g v_n=-\rho_{2,n}\left(\frac{e^{-u_n}}{\int_{\Sigma}e^{-u_{n}}dV_g}-1\right),\;\;\;\bar v_n=0
\end{equation}
satisfies\;$v_n$\;is uniformly bounded in \;$W^{2,p}$\; (thanks to standard elliptic regularity). Thus taking \;$p$\;so large we have by Sobolev-Embedding theorem \;$v_n$\;is bounded in \;$C^{1,\alpha}$. Now defining \;$w_n$\;by \;$w_{n}=u_n-v_n$, we have that \;$w_n$\;sove the folowing PDE
\begin{equation}
-\D_g w_n=\rho_{1,n}\left(e^{v_n}e^{v_n}\frac{e^{u_n}}{\int_{\Sigma}e^{v_n}e^{u_{n}}dV_g}-1\right),\;\;\;\bar w_n=0
\end{equation}
So using Theorem\;$\ref{th:yyli}$, we get 
\;$w_n$\;is uniformly bounded in \;$L^{\infty}(\Sig)$. Thus we get \;$u_n$\;is uniformly bounded in \;$L^{\infty}(\Sig)$. Hence, by standard elliptic regularity theory we get \;$u_n$\;is bounded in $C^m(\Sigma)$\;for every positive integer \;$m$. Hence the proposition is proved.
\end{pf}

\section{Proof of Theorem\;$\ref{eq:theo}$}
This section deals with the proof of Theorem\;$\ref{eq:theo}$. It is divided into four subsections. The first one is concerned with the definition of the {\em formal barycenters} of \;$\Sig$, and some related results. The second one is about the derivation of an improvemment of the Moser-Trudinger type inequality given by Theorem\;$\ref{th:os1}$\;and corollaries. The third one deals with the construction of a continuous map from large negative sublevels of \;$\text{II}_{\rho}$ \;into \;$\Sig_{k}$ (for the definition see Subsection 1)\;and an other one from\;$\Sig_k$\; into suitable negative sublevels of \;$\text{II}_{\rho}$. The last one describes the topological argument.
\subsection{Barycenters and Properties}
As said in the introduction of the Section, we start by recalling the definition of the so called {\em formal barycenters} of \;$\Sig $.\\ For $k \in \N$,
we let $\Sig_k$ denote the family of formal sums
\begin{equation}\label{eq:Mk}
    \Sig_k = \sum_{i=1}^k t_i \d_{x_i}; \qquad \qquad t_i \geq 0,
    \quad \sum_{i=1}^k t_i = 1, \quad x_i \in \Sig,
\end{equation}
where $\d_x$ stands for the Dirac delta at the point $x \in \Sig$.
We endow this set with the weak topology of distributions. This is
known in literature as the {\em formal set of barycenters} of $\Sig$
(of order $k$), see \cite{bah}, \cite{bc}, \cite{bred}. Although
this is not in general a smooth manifold (except for $k = 1$), it is
a {\em stratified set}, namely union of cells of different
dimensions. The maximal dimension is $3k - 1$, when all the points
$x_i$ are distinct and all the $t_i$'s belong to the open interval
$(0,1)$.
\vspace{6pt}

\noindent
After introducing the set of formal barycenters, we give the following well-know result (see \cite{dm}) which is necessary for the topological argument below.
\begin{lem}\label{l:nonco} (well-known)
For any $k \geq 1$ one has $H_{3k-1}(\Sig_k;\Z_2) \neq 0$. As a
consequence $\Sig_k$ is non-contractible.
\end{lem}
Next we introduce a distance on\;$\Sig_k$.\\
If $\var \in C^1(\Sig)$ and if $\s \in \Sig_k$, we denote
the action of $\s$ on $\var$ as
$$
  \langle \s, \var \rangle = \sum_{i=1}^k t_i \var(x_i), \qquad \quad
\s = \sum_{i=1}^k t_i \d_{x_i}.
$$
Moreover, if $f$ is a non-negative $L^1$ function on $\Sig$ with
$\int_\Sig f dV_g = 1$, we can define a distance of $f$ from
$\Sig_k$ in the following way
\begin{equation}\label{eq:distf}
    dist(f, \Sig_k) = \inf_{\s \in \Sig_k} \sup \left\{ \left|
  \int_{\Sig} f \var dV_g - \langle \s, \var \rangle \right|
  \; | \; \|\var\|_{C^1(\Sig)} = 1 \right\}.
\end{equation}
We also let
$$
  \mathcal{D}_{\e,k} = \left\{ f \in L^1(\Sig) \; : \; f \geq 0,
  \|f\|_{L^1(\Sig)} = 1, dist(f, \Sig_k) < \e \right\}.
$$

\noindent From a straightforward adaptation of the arguments of
Proposition 3.1 in \cite{dm}, we obtain the following result.

\begin{pro}\label{p:projdm} Let $k$ be a positive integer, and for
$\e > 0$ let $\mathcal{D}_{\e,k}$ be as above. Then there exists
$\e_k > 0$, depending on $k$ and $\Sig$ such that, for $\e \leq
\e_k$ there exists a continuous map $\Pi_k : \mathcal{D}_{\e,k} \to
\Sig_k$.
\end{pro}

\subsection{Improved Moser-Trudinger inequality and applications}
In this subsection we analyze the Moser-Trudinger type inequality given by Theorem\;$\ref{th:os1}$. We prove that depending on the amount of concentration of \;$e^{u}$ it get an improvement. From this we charcterizes low sublevels of \;$\text{II}_{\rho}$\;in terms of the concentration of \;$e^{u}$.

\begin{pro}\label{p:imprc}
Let $\d_0 > 0$, $\ell \in \N$, and let $S_{1}, \dots, S_{\ell}$ be
subsets of $\Sig$ satisfying $dist(S_{i},S_{j}) \geq \d _0$ for $i
\neq j$. Let $\g_0 \in \left( 0, \frac{1}{\ell} \right)$. Then, for
any $\tilde{\e} > 0$ there exists a constant $C = C(\tilde{\e},
\d_0, \g_0, \ell, \Sig)$ such that
$$
  \ell \log \int_\Sig e^{(u - \ov{u})} dV_g +
  \log \int_\Sig e^{-(u - \ov{u})} dV_g \leq C + \frac{1}{16 \pi -
  \tilde{\e}}\int_{\Sig} |\n u|^2dV_g
$$
provided the function $u$ satisfies  the relations
\begin{equation}\label{eq:ddmm}
    \frac{\int_{S_{i}} e^{u} dV_g}{\int_\Sig e^{u} dV_g}
    \geq \g_0, \quad \qquad i \in \{1, \dots, \ell\}.
\end{equation}
\end{pro}
Before making the proof we recall the following Lemma whose proof is a trivial adaptation of Lemma 3.2 in \cite{mn}.

\begin{lem}\label{l:step1}
Under the assumptions of Proposition \ref{p:imprc}, there exist
numbers $\tilde{\g}_0, \tilde{\d}_0 > 0$, depending only on $\g_0,
\d_0, \Sig$, and $\ell$ sets $\tilde{S}_1, \dots, \tilde{S}_{\ell}$
such that $d(\tilde{S}_i, \tilde{S}_j) \geq \tilde{\d}_0$ for $i
\neq j$ and such that
$$
  \frac{\int_{\tilde{S}_1} e^{u} dV_g}{\int_\Sig e^{u} dV_g} \geq \tilde{\g}_0,
  \quad \frac{\int_{\tilde{S}_1} e^{u} dV_g}{\int_\Sig e^{u} dV_g}
    \geq \tilde{\g}_0; \qquad  \qquad \quad
  \frac{\int_{\tilde{S}_i} e^{-u} dV_g}{\int_\Sig e^{-u} dV_g}
    \geq \tilde{\g}_0, \quad i \in \{2, \dots, \ell\}.
$$
\end{lem}

\begin{pfn} {\sc of Proposition \ref{p:imprc}.} We use the argument in \cite{dm} adapted to our purpose.
. Firts of all let $\tilde{S}_1, \dots,
\tilde{S}_{\ell}$ be given by Lemma \ref{l:step1}. Moreover without
loss of generality we assume that $\ov{u} = 0$.  We have there exist
$\ell$ functions $g_1, \dots, g_{\ell}$ satisfying the properties
\begin{equation}\label{eq:gi}
    \left\{%
\begin{array}{ll}
    g_i(x) \in [0,1] & \hbox{ for every } x \in \Sig; \\
    g_i(x) = 1, & \hbox{ for every } x \in \tilde{S}_i, i = 1, \dots, \ell; \\
    supp(g_i) \cap supp(g_j) = \emptyset, & \hbox{ for } i \neq j; \\
    \|g_i\|_{C^2(\Sig)} \leq C_{\tilde{\d}_0}, &  \\
\end{array}%
\right.
\end{equation}
where $C_{\tilde{\d}_0}$ is a positive constant depending only on
$\tilde{\d}_0$.\\
Next we decompose the function $u$  in Fourier mode (to be choosen later) as follows 
\begin{equation}\label{eq:decui}
    u = \hat{u} + \tilde{u};  \qquad \qquad \hat{u}\in
  L^\infty(\Sig).
\end{equation} 
Now using Lemma \ref{l:step1}, for any
$b \in 2, \dots, \ell$ we can write that
\begin{eqnarray*}
  \ell \log \int_\Sig e^{u} dV_g + \log \int_\Sig
  e^{-u} dV_g & = & \log \left[ \left( \int_\Sig e^{u} dV_g
   \int_\Sig e^{-u} dV_g \right) \left( \int_\Sig
   e^{u} dV_g \right)^{\ell - 1} \right] \\
   & \leq & \left[ \left( \int_{\tilde{S}_1} e^{u} dV_g
   \int_{\tilde{S}_1} e^{-u} dV_g \right) \left(
   \int_{\tilde{S}_b} e^{u} dV_g \right)^{\ell - 1}
   \right] - \ell \log \tilde{\g}_0 \\
   & \leq & \log \left[ \left( \int_\Sig e^{g_1 u} dV_g
   \int_\Sig e^{-g_1 u} dV_g \right) \left( \int_\Sig
   e^{g_b u} dV_g \right)^{\ell - 1} \right] \\ & - &
   \ell \log \tilde{\g}_0,
\end{eqnarray*}
Using the fact that $\hat{u}$  belong to
$L^\infty(\Sig)$, wearrive to 
\begin{eqnarray*}
  \ell \log \int_\Sig e^{u} dV_g + \log \int_\Sig
  e^{-u} dV_g & \leq & \log \left[ \left( \int_\Sig e^{g_1 \tilde{u}} dV_g
   \int_\Sig e^{-g_1 \tilde{u}} dV_g \right) \left( \int_\Sig
   e^{g_b \tilde{u}} dV_g \right)^{\ell - 1} \right] \\ & - &
   \ell \log \tilde{\g}_0 + (\ell+1) \|\hat{u}\|_{L^\infty(\Sig)}
   .
\end{eqnarray*}
Thus we get
\begin{eqnarray}\label{eq:aa} \nonumber
  \ell \log \int_\Sig e^{u} dV_g + \log \int_\Sig e^{u} dV_g
   & \leq & \log \int_\Sig e^{g_1 \tilde{u}} dV_g + \log
   \int_\Sig e^{-g_1 \tilde{u}} dV_g + (\ell - 1) \int_\Sig e^{g_b \tilde{u}}
   dV_g \\ & - & \ell \log \tilde{\g}_0 + (1+\ell )\|\hat{u}\|_{L^\infty(\Sig)} .
\end{eqnarray}
Now apply Theorem\;$\ref{th:os1}$\; with parameters $(8\pi,8\pi)$\ to the couple $(g_1 \tilde{u}, -g_1
\tilde{u})$, and the standard Moser-Trudinger inequality
$\eqref{eq:mt}$\; to\; $g_b \tilde{u}$ \; we obtain 
\begin{eqnarray}\label{eq:bb} \nonumber
   \log \int_\Sig e^{g_1 \tilde{u}} dV_g + \log \int_\Sig
  e^{-g_1 \tilde{u}} dV_g & \leq & \frac{1}{16\pi} \int_{\Sig}|\n g_1u|^2dV_g +C;
\end{eqnarray}
$$
  (\ell - 1) \int_\Sig e^{g_b \tilde{u}} dV_g \leq \frac{(\ell -
  1)}{16 \pi} \int_\Sig |\n (g_b \tilde{u})|^2 dV_g + (\ell - 1)
  \ov{g_b \tilde{u}} + (\ell - 1) C.
$$
Putting together \eqref{eq:aa}-\eqref{eq:bb} we get
\begin{eqnarray}\label{eq:dd} \nonumber
  \ell \log \int_\Sig e^{u} dV_g + \log \int_\Sig e^{-u} dV_g\leq \frac{1}{16\pi} \int_{\Sig}|\n g_1u|^2dV_g +\frac{(\ell -
  1)}{16 \pi} \int_\Sig |\n (g_b \tilde{u})|^2 dV_g + (\ell - 1)
  \ov{g_b \tilde{u}} + C.
   \nonumber
\end{eqnarray}
Next, by  interpolation, for any $\e > 0$ there exists
$C_{\e,\tilde{\d}_0}$ (depending only on $\e$ and $\tilde{\d}_0$)
such that
\begin{equation} 
\frac{1}{16\pi}\int_{\Sig}|\n (g_i\tilde {u})|^2dV_g \leq \frac{1}{16\pi}\int_{\Sig}g_i^2|\n \tilde {u}|^2dV_g
  + \frac{\epsilon}{16\pi}\int_{\Sig}|\n \tilde {u}|^2dV_g+C_{\e,\tilde{\d}_0} \int_{\Sig} \tilde{u}^2  dV_g.
\end{equation}
Hence inserting this inequality into \;$\eqref{eq:dd}$\; we get
\begin{eqnarray*}
\ell \log \int_\Sig e^{u} dV_g + \log \int_\Sig e^{-u} dV_g\leq \frac{1}{16\pi}\int_{\Sig}g_1^2|\n \tilde {u}|^2dV_g 
   +\frac{(\ell -
  1)}{16 \pi} \int_{\Sig} g_b^2|\n  \tilde{u}|^2 dV_g + \frac{\ell \epsilon}{16\pi}\int_{\Sig}|\n \tilde{u}|^2dV_g\\+\ell C_{\e,\tilde{\d}_0} \int_{\Sig} \tilde{u}^2  dV_g +
  (\ell-1)
  \ov{g_b \tilde{u}} + C
  ,
\end{eqnarray*}
Now for $b = 2, \dots, \ell$, we  choose $b \in \{ 2, \dots, \ell \}$ such that
$$
   \frac{1}{16\pi}\int_{\Sig}g_b^2|\n \tilde {u}|^2dV_g\leq \frac{1}{\ell - 1} \frac{1}{16\pi} 
  \int_{\cup_{s=2}^{\ell} supp(g_s)}  |\n u|^2dV_g.
$$
On the other hand since  the $g_i's$ have disjoint supports, see \eqref{eq:gi}, then
last formula yields
\begin{eqnarray*}
\ell \log \int_\Sig e^{u} dV_g + \log \int_\Sig e^{-u} dV_g\leq \frac{1+\ell \epsilon}{16\pi}\int_{\Sig}|\n \tilde {u}|^2dV_g 
   +\ell C_{\e,\tilde{\d}_0} \int_{\Sig} \tilde{u}^2  dV_g +
  (\ell-1)
  \ov{g_b \tilde{u}} + C
\end{eqnarray*}
Next, by elementary estimates we find
\begin{eqnarray*}
  \ell \log \int_\Sig e^{u} dV_g + \log \int_\Sig e^{-u}
  dV_g \leq  \frac{1+\ell \epsilon}{16\pi}\int_{\Sig}|\n \tilde {u}|^2dV_g +
  C_{\e, \tilde{\d}_0, \ell} \int_{\Sig} \tilde{u}^2 dV_g
  +C_{\e, \tilde{\d}_0, \ell,\tilde{\g}_0}
   + \ell \|\hat{u}\|_{L^\infty(\Sig)} .
\end{eqnarray*}
Now comes the choice of $\hat{u}$, see
\eqref{eq:decui}. We choose $\tilde{C}_{\e, \tilde{\d}_0, \ell}$ to
be so large that the following property holds
$$
  C_{\e, \tilde{\d}_0, \ell} \int_{\Sig} v^2  dV_g <
\frac{\epsilon}{16\pi}\int_{\Sig}|\n v|^2dV_g, \qquad \forall
v \in V_{\e, \tilde{\d}_0, \ell},
$$
where $V_{\e, \tilde{\d}_0, \ell}$ denotes the span of the
eigenfunctions of the Laplacian on $\Sig$ corresponding to
eigenvalues bigger than $\tilde{C}_{\e, \tilde{\d}_0, \ell}$.

Then we set
$$
  \tilde{u} = P_{V_{\e, \tilde{\d}_0, \ell}} u; \qquad \qquad
  \hat{u} = P_{V_{\e, \tilde{\d}_0, \ell}^\perp} u,
$$
where $P_{V_{\e, \tilde{\d}_0, \ell}}$ (resp. $P_{V_{\e,
\tilde{\d}_0, \ell}^\perp}$) stands for the orthogonal projection
onto $V_{\e, \tilde{\d}_0, \ell}$ (resp. $V_{\e, \tilde{\d}_0,
\ell}^\perp$). Since $\ov{u} = 0$, the $H^1$-norm and the
$L^\infty$-norm on $V_{\e, \tilde{\d}_0, \ell}^\perp$ are equivalent (with
a proportionality factor which depends on $\e, \tilde{\d}_0$ and
$\ell$), hence by our choice of $u$  there holds
$$
  \|\hat{u}\|_{L^\infty(\Sig)}^2 \leq \hat{C}_{\e, \tilde{\d}_0, \ell}
  \frac{1}{16\pi} \int_{\Sig}|\n \hat u|^2dV_g dV_g; \qquad \quad C_{\e, \tilde{\d}_0, \ell} \int_{\Sig}
\tilde{u}^2 dV_g  < \frac{\e}{16\pi} \int_{\Sig}|\n \tilde u|^2dV_g.
$$
Hence the last formulas imply
\begin{eqnarray*}
  \ell \log \int_\Sig e^{u} dV_g + \log \int_\Sig e^{-u}
  dV_g & \leq & \frac{1}{16\pi} (1 + 3 \ell \e)\int_{\Sig}|\n u|^2dV_g + \hat{C}_{\e, \tilde{\d}_0, \ell,\tilde{\g}_0}.
\end{eqnarray*}
This concludes the proof.
\end{pfn}
\vspace{6pt}

\indent
In the remaining of this subsection we will apply the above Proposition to understand the structure of the sublevels of \;$\text{II}_{\rho}$. Before this we state a Lemma  which gives sufficient conditions for the improvement to hold. Its proof can be found in \cite{dm}.
\begin{lem}\label{l:er}
Let $f \in L^1(\Sig)$ be a non-negative function with
$\|f\|_{L^1(\Sig)} = 1$. We also fix an integer $\ell$ and suppose
that the following property holds true. There exist $\e
> 0$ and $r > 0$ such that
$$
  \int_{\cup_{i=1}^\ell B_r(p_i)} f dV_g < 1 - \e \qquad \qquad \hbox{ for all
the $\ell$-tuples }
  p_1, \dots, p_\ell \in \Sig.
$$
Then there exist $\ov{\e} > 0$ and $\ov{r} > 0$, depending only on
$\e, r, \ell$ and $\Sig$ (and not on $f$), and $\ell + 1$ points
$\ov{p}_1, \dots, \ov{p}_{\ell+1} \in \Sig$ (which depend on $f$)
satisfying
$$
  \int_{B_{\ov{r}}(\ov{p}_1)} f dV_g > \ov{\e}, \; \dots, \;
  \int_{B_{\ov{r}}(\ov{p}_{\ell+1})} f dV_g > \ov{\e}; \qquad \qquad
  B_{2 \ov{r}}(\ov{p}_i) \cap B_{2 \ov{r}}(\ov{p}_j) = \emptyset
  \hbox{ for } i \neq j.
$$
\end{lem}

\begin{pro}\label{l:II<-M}
Suppose $\rho_1 \in (8 \pi k, 8 \pi (k + 1))$ and that $\rho_2 < 8
\pi$. Then for any $\e > 0$ and any $r
> 0$ there exists a large positive $L = L(\e, r)$ such that for every
$u \in H^1(\Sig) $ with $\text{II}_{\rho}(u) \leq -
L$ and with $\int_\Sig e^{u} dV_g = 1$,  there exists
$k$ points $p_{1,u}, \dots, p_{k,u} \in \Sig$ such that
\begin{equation}\label{eq:caz}
 \int_{\Sig \setminus \cup_{i=1}^{k} B_r(p_{i,u})} e^{u} dV_g <
 \e.
\end{equation}
\end{pro}
\begin{pf}
To prove the proposition. we willl argue by contradiction. So suppose it does not holds, then aplying Lemma\;$\ref{l:er}$\; with \;$l=k$\; and\;$f=e^u$\;, we have that there exists \;$\delta_0$, \;$\gamma_0$\; and set \;$S_1,\cdots, S_{l+1}$\;such that  \;$d(S_i,S_j)\geq \delta_0$ and 
\begin{equation}
\frac{\int_{S_i}e^udV_g}{\int_{\Sig}e^udV_g}\geq \gamma_0,\;\;\text{for}\;\;i=1,\cdots,l+1.
\end{equation}
Next from Jensen's inequality and the fact that \;$\int_{\Sig}e^udVg=1$\;we get 
\begin{equation}\label{eq:imp}
\bar u\leq 0\;\;\text{and}\;\;\;\log\int_{\Sig}e^{-u+\bar u}dV_g\geq 0.
\end{equation}
Now since \;$\rho_1<8\pi(k+1)$\;and \;$\rho_2<4\pi$\; then there exits a small \;$\tilde \epsilon >0$\;such that 
\begin{equation}\label{eq:imp2}
(16\pi-\tilde \epsilon)(k+1)>2\rho_1\;\;\;\text{and}\;\;16\pi-\tilde \epsilon>2\rho_2.
\end{equation}
On the other hand  from the definition of \;$\text{II}_{\rho}$\;we have that
\begin{equation}
\begin{split}
\text{II}_{\rho}(u)=\frac{1}{2}\int_{\Sig}|\n u|^2dV_g-(\frac{16\pi-\tilde \epsilon}{2})(k+1)\log\int_{\Sig}e^{u-\bar u}dV_g-(\frac{16\pi-\tilde \epsilon}{2})\log\int_{\Sig}e^{-u+\bar u}dV_g\\+(\frac{16\pi-\tilde \epsilon}{2})(k+1)-\rho_1)\log\int_{\Sig}e^{u-\bar u}dV_g+(\frac{16\pi-\tilde \epsilon}{2}-\rho_2)\log\int_{\Sig}e^{-u+\bar u}dV_g.
\end{split}
\end{equation}

Hence using\;$\eqref{eq:imp}$, the normalizatiuon \;$\int_{\Sig}e^udV_g=1$\; and \;$\eqref{eq:imp2}$, we  get
\begin{equation}
\text{II}_{\rho}(u)\geq \frac{1}{2}\int_{\Sig}|\n u|^2dV_g-(\frac{16\pi-\tilde \epsilon}{2})(k+1)\log\int_{\Sig}e^{u-\bar u}dV_g-(\frac{16\pi-\tilde \epsilon}{2})\log\int_{\Sig}e^{-u+\bar u}dV_g.
\end{equation}
Next using \;$\eqref{p:imprc}$\;we obtain
\begin{equation}
\text{II}_{\rho}(u)\geq -C
\end{equation}
Hence the proposition is proved.
\end{pf}
\vspace{6pt}

\indent
The next result is a direct corollary of  Proposition\;$\ref{l:II<-M}$. It gives the distance of \;$e^u$\; from \;$\Sig_k$\; for \;$u$\;belonging to low sublels of \;$\text{II}_{\rho}$\; and \;$\int_{\Sig}e^udV_g=1$
\begin{cor}\label{eq:impmt}
Let $\ov{\e}$ be a (small) arbitrary positive number and \;$k$\;be given as in Theorem\;$\ref{eq:theo}$. Then there
exists $L > 0$ such that, if \;$\text{II}(u) \leq - L$ \;and \;
$\int_\Sig e^{u} dV_g = 1$, then  we have that\;$d(e^{u}, \Sig_k) \leq 
\ov{\e}$.
\end{cor}
\begin{pf}
Let \;$\epsilon>0$, \;$r>0$\;(to be fixed later) and let \;$L$\;be the corresponding constant given by Proposition\;$\ref{l:II<-M}$.  We let \;$p_{1},\cdots,p_{k}$\;be the points given by Proposition\;$\ref{l:II<-M}$\; and we define \;$\s\in \Sig_k$\;as follows
\begin{equation}
\s=\sum_{i=1}^{k}t_i\delta_{p_{i}}\;\;\text{where}\;t_i=\int_{A_{r,i}}e^{4u}dV_j,\;\;A_{r,i}:=B_{p_i}(r)\setminus\cup_{s=1}^{i-1}B_{p_s}(r),\;i=1,\cdots,k-1,\; t_k=1-\sum_{i=1}^{k-1}t_i.
\end{equation}
By construction we have \;$A_{r,i}$\;are disjoint and \;$\cup_{i=1}^{k-1} A_{r,i}=\cup_{i=1}^{k-1} B_{p_i}(r)$. Now let \;$\varphi\in C^1(\Sig)$\;be such that \;$||\varphi||_{C^1(\Sig)}=1$,\; we have that by triangle inequality
\begin{equation}
\left|\int_{\Sig}e^u\varphi-<\s,\varphi>\right|\leq\sum_{i=1}^{k-1}\left|\int_{A_{r,i}}e^{u}(\varphi-\varphi(p_i)\right|+\left|\int_{\Sig\setminus \cup_{i=1}^{k-1}A_{r,i}}e^{u}(\varphi-,\varphi(p_k)\right|.
\end{equation}
Thus by using Mean value formula and \;$\eqref{eq:caz}$\;we get
\begin{equation}
\left|\int_{\Sig}e^{4u}\varphi-<\s,\varphi>\right|\leq C_{\Sig}r+C_{\Sig}r\epsilon.
\end{equation}
So by choosing \;$\epsilon$\;and \;$r$\;so small that\;$ C_{\Sig}r+C_{\Sig}r\epsilon <\bar\epsilon$, and recalling that \;$d$\;is the metric given by\;$C^1(\Sig)^*$, we obtain
\begin{equation}
d(e^{u},\Sig_k)<\bar \epsilon;
\end{equation}
hence we are done.
\end{pf}

\subsection{Construction of the projections \;$\Psi$\;and \;$\Phi$}
In this Subsection we construct two global continuous non-trivial projections in order to show that large negative sublevels of \;$\text{II}_{\rho}$\;have the same homology as \;$\Sig_k$, see Proposition\;$\ref{eq:test}$\;below.

\begin{pro}\label{p:glob}
Let \;$k$,\;$\rho_1$\;and \;$\rho_2$\;as in Theorem\;$\eqref{eq:theo}$. Then there exists
a large $L > 0$ and a continuous projection $\Psi$ from $\{ \text{II}_\rho
\leq - L \} \cap \left\{ \int_\Sig e^{u} dV_g = 1 \right\}$ (with
the natural topology of $H^1(\Sig)$ ) onto
$\Sig_{k}$ which is homotopically non-trivial.
\end{pro}

\begin{pf}
We fix $\e_{k}$ so small that Proposition \ref{p:projdm} applies
. Then we apply Corollary \ref{eq:impmt} with $\ov{\e} =
\e_k$. We let $L$ be the corresponding large number, so that if
\;$\text{II}_{\rho}(u) \leq - L$, then $dist(e^{u}, \Sig_{k}) < \e_k$. Hence
for these ranges of $u$ , since the map $u \mapsto e^u$
is continuous from $H^1(\Sig)$ into $L^1(\Sig)$, setting \;$\Psi(u)=\Pi_k(e^u)$ (where\;$\Pi_k$\;is given by Proposition \ref{p:projdm}), we have \;$\Psi(\cdot)$\;is continuous. The non-triviality of this map is a consequence of
Proposition \ref{eq:test} {\bf (ii)}.
\end{pf}
\vspace{6pt}

\indent
Next, we show that one can map \;$\Sig_k$\;into very large negative sublevels of \;$\text{II}_{\rho}$. To do this we start by introducing some notations..\\
Given \;$\s=\sum_{i=1}^kt_i\d_{x_i}\in \Sig_k$\; and \;$\l$\;a positive real number, we set
\begin{equation}\label{eq:fi}
\varphi_{\s,\l}(y)=\log\sum_{1=1}^k\left(\frac{\l}{1+\l^2d_i(y)^2}\right)^2-\log\pi,\;\;\;y\in \Sig;
\end{equation}
where\;$d_i(y)=d(y,x_i)$.\\
We remark that, since the distance function is lipschitz, then \;$\varphi_{\s,\l}$\;is, hence due to Sobolev embedding is an element of \;$H^1(\Sig)$.\\
We have the following Proposition about \;$\varphi_{\s,\l}$.
\begin{pro}\label{eq:test}
Supposs \;$k$,\;$\rho_1$\;and \;$\rho_2$\;as in \;Theorem\;$\ref{eq:theo}$. For \;$\l>0$\;and \;$\s\in \Sig_k$\;we define\;$$\Phi_{\l}:\Sig_k\rightarrow H^1(\Sig)$$\;as 
\begin{equation*}
\Phi_{\l}(\s)=\varphi_{\s,\l}
\end{equation*}
where \;$\varphi_{\s,\l}$\;is as in\;$\eqref{eq:fi}$. Then for \;$L$\;suficciently the exist \;$\bar \l>0$\;such that

$
\noindent (i)\; \text{II}_{\rho}(\Phi_{\l}(\s))\leq -L\;\;\text{uniformly in }\;\;\s\in \Sig_k\;\;\;\l\geq \bar \l;
$
\\\\
$
\noindent (ii)\; \Psi\circ \Phi _{\l}\;\;\text{is homotopic to the identity on}\;\;\Sig_k\;\;\text{for }\;\; \;\l\;\;\text{large}..
$
\end{pro}
\begin{pf}
To prove \;$(i)$,we first claim that as \;$\l\rightarrow +\infty$\; the following estimate holds
\begin{equation}\label{eq:mean}
\int_{\Sig}\varphi_{\s,\l}=-2(1+o_{\l}(1))\log \l,
\end{equation}
\begin{equation}\label{eq:expo}
\log\int_{\Sig}e^{\varphi_{\s,\l}}dV_g=O(1)\;\;\text{and}\;\;\log\int_{\Sig}e^{-\varphi_{\s,\l}}dV_g=2(1+o_{\l}(1))\log\l,
\end{equation}
and
\begin{equation}\label{eq:grad}
\int_{\Sig}|\n \varphi_{\s,\l}|^2dV_g\leq 32k\pi(1+o_{\l}(1))\log \l;.
\end{equation}
{\bf Proof of Claim}\\
{\em Proof of \;$\eqref{eq:mean}$}\\
 Let \;$\delta\in (0,diam(\Sig))$\; be small. We have that
 \begin{equation}
 2\log\frac{\l}{1+\l^2diam(\Sig)^2}-\log\pi\leq \varphi_{\s,\l}\leq 2\log\frac{\l}{1+\l^2\delta^2}-\log \pi\;\;\;\text{for}\;\;y\in \Sig\setminus \cup_{1=1}^kB_{x_i}(2\delta);
 \end{equation}
 and
 \begin{equation}\label{eq:rewrite}
  2\log\frac{\l}{1+4\l^2\delta^2}-\log \pi\leq \varphi_{\s,\l} \leq 2\log\l-\log \pi\;\;\;\text{for}\;\;y\in \cup_{1=1}^kB_{x_i}(2\delta);
 \end{equation}
Now rewritting \;$\eqref{eq:rewirte}$\;we obtain
\begin{equation*}
-2\log\l-2\log\left(1+\frac{diam(\Sig)^2}{\l^2}\right)-\log\pi\leq \varphi_{\s,\l}\leq -2\log\l-2\log\left(1+\frac{\delta^2}{\l^2}\right)-\log\pi\;\;\in \Sig\setminus \cup_{1=1}^kB_{x_i}(2\delta)
\end{equation*}
Thus combining all, get
\begin{equation}
\int_{\Sig}\varphi_{\s,\l}dV_g=-2\log\l(1+O(\delta^2))+O(1)+O(\delta^2)(|\log\l|+|\log\delta|)
\end{equation}
Hence letting \;$\delta$\;tends to zero we get the desired conclusion.\\\\
{\em Proof of \;$\eqref{eq:expo}$}\\
The proof of $(ii)$ \;comes from direct calculations.\\.\\\\
{\em Proof of \;$\eqref{eq:grad}$}

The proof of this inequality relies on showing the following two pointwise
estimates on the gradient of $\var_{\l,\s}$
\be
  |\n \var_{\l,\s}(y)| \leq C \l; \qquad \qquad \hbox{ for every } y
\in \Sig, \label{eq:ptwbd} \ee where $C$ is a constant independent
of $\s$ and $\l$, and \be \label{eq:dbd} |\n \var_{\l,\s}(y)| \leq
\frac{4}{d_{min}(y)} \qquad \hbox{ where } \qquad  d_{min}(y) =
\min_{i=1,\dots,m} d(y,x_i). \ee

For proving \eqref{eq:ptwbd} we notice that the following inequality
holds
\begin{equation}\label{eq:ineq}
    \frac{\l^2 d(y,x_i)}{1 + \l^2 d^2(y,x_i)} \leq C \l, \qquad i =
1, \dots, m,
\end{equation}
where $C$ is a fixed constant (independent of $\l$ and $x_i$).
Moreover we have
\begin{equation}\label{eq:grad1}
    \n \var_{\l,\s}(y) = - 2 \l^2 \frac{\sum_i t_i (1 + \l^2 d_i^2(y))^{-3}
\n_y (d_i^2(y))}{\sum_j t_j (1 + \l^2 d_j^2(y))^{-2}}.
\end{equation}
Using the fact that $|\n_y (d_i^2(y))| \leq 2 d_i(y)$ and inserting
\eqref{eq:ineq} into \eqref{eq:grad1} we obtain immediately
\eqref{eq:ptwbd}. Similarly we find
\begin{eqnarray*}
  |\n \var_{\l,\s}(y)| & \leq & 4 \l^2 \frac{\sum_i t_i (1 + \l^2 d_i^2(y))^{-3}
d_i(y)}{\sum_j t_j (1 + \l^2 d_j^2(y))^{-2}} \leq 4 \l^2
\frac{\sum_i t_i (1 + \l^2 d_i^2(y))^{-2}
\frac{d_i(y)}{\l^2 d_i^2(y)}}{\sum_j t_j (1 + \l^2 d_j^2(y))^{-2}} \\
   & \leq & 4 \frac{\sum_i t_i (1 + \l^2 d_i^2(y))^{-2}
\frac{1}{d_{min}(y)}}{\sum_j t_j (1 + \l^2 d_j^2(y))^{-2}} \leq
\frac{4}{d_{min}(y)},
\end{eqnarray*}
which is \eqref{eq:dbd}.\\
From we infer that 
\begin{equation}\label{eq:grad1}
\int_{\cup_{i=1}^kB_{x_i}(\frac{1}{\l})}|\n_g \varphi_{\s,\l}|^2\leq Ck;
\end{equation}
for some constant depending only on \;$\Sig$.\\
Now for every \;$i=1,\cdots,k$\;we set
\begin{equation}
B_i=\{y\in \Sig \;d(y,x_i)=d_{min}(y)\}
\end{equation}
and we have
\begin{equation}
\begin{split}
\int_{\Sig\cup_{i=1}^kB_{x_i}(\frac{1}{\l})}|\n_g \varphi_{\s,\l}|^2dV_g\leq \sum_{i=1}^k\int_{B_i\setminus B_{x_i}(\frac{1}{\l})}|\n_g \varphi_{\s,\l}|^2dV_g\leq 16\sum_{i=1}^k\int_{B_i\setminus B_{x_i}(\frac{1}{\l})}\frac{1}{d(y,x_i)^2}dVg(y)\\
\leq 32\pi(1+o_{\l}(1))\log\l+O(1).
\end{split}
\end{equation}
From this and \;$\eqref{eq:grad1}$\;we deduce\;$\eqref{eq:grad}$.\\\\
Hence the proof of Claim is complete.\\
Next using the Claim and the definition of \;$\text{II}_{\rho}$\;we get 
\begin{equation}
\text{II}_{\rho}(\varphi_{\s,\l}) \leq (16k\pi-2\rho_1+o_{\l}(1)\log\l+O(1).
\end{equation}
Thus using the fact that \;$8k\pi<\rho_1$\; we get that
\begin{equation}
\text{II}_{\rho}(\varphi_{\s,\l})\rightarrow -\infty \;\;\text{uniformly in}\;\s
\end{equation}
Hence the proof of\;$(i)$\;is completed. \\
Now let us show\;$(ii)$. Firts of all we remark for every given \;$x$, the trivial convergence holds
\begin{equation}
(\frac{\l}{1+\l^2d(x,y)^2})^2\rightharpoonup \pi\d_x
\end{equation}
in the weak sens of measure. Hence using the definition of \;$\varphi_{\s,\l}$\; one check easily that
\begin{equation}\label{eq:cdelta}
e^{\varphi_{\s,\l}}\rightharpoonup \s.
\end{equation}
Onb the other hand from \;$(i)$\;we have that the following composition  for large \;$\l$
\begin{equation}
T_{\l}=\Psi\circ\Phi_{\l}
\end{equation}
is well defined. Moreover from $\eqref{eq:cdelta}$ and the continuity of \;$\Psi$\;we infer that for \;$\bar \l$\;large \;$T_{\l}$\;is an homotopy beetween \;$\Psi\circ\Phi_{\bar \l}$\;and identity on \;$\Sig_k$. Thus the proof of \;$(ii)$\;is complete. Hence the proof of the proposition is concluded.

\end{pf}
\subsection{Topological argument}
In this Subsetion we perform the topological argument in order to produce solutions. We will employ a min-max scheme based on the topological cone \;$C_k$ (for precise definition see  below)\;over\;$\Sig_k$. As anticipated in the introduction, we then
define a modified functional $\text{II}_{t \rho_1, t \rho_2}$ for which we
can prove existence of solutions in a dense set of the values of
$t$. Following an idea of Struwe, this is done proving the a.e.
differentiability of the map $t \mapsto \a_{t \rho}$, where $\a_{t
\rho}$ is the minimax value for the functional $\text{II}_{t \rho_1, t
\rho_2}$ given by the scheme.\\
Let \;$C_k$\;be the topological cone over \;$C_k$, see. First, let $L$ be so large that Proposition
\ref{p:glob} applies with $\frac L4$, and choose then $\Phi$ such
that Proposition \ref{eq:test} applies for $L$. Fixing $L$ and
$\Phi$, we define the class of maps
\begin{equation}\label{eq:PiPi}
    \Pi_{\Phi_{\bar \l}} = \left\{ \pi : C_{k} \to H^1(\Sig) 
 \; : \; \pi \hbox{ is continuous and } \pi|_{\Sig_{k} (=
\partial K_{k})} = \Phi_{\bar \l} \right\}.
\end{equation}

Then we have the following properties.

\begin{lem}\label{l:minmax}
The set $\Pi_{\Phi}$ is non-empty and moreover, letting
$$
  \a_\rho = \inf_{\pi \in \Pi_{\Phi}}
  \; \sup_{m \in C_k} \text{II}_{\rho_1,\rho_2}(\pi(m)), \qquad
  \hbox{ there holds } \qquad \a_\rho > - \frac
  L2.
$$
\end{lem}

\begin{pf}
To prove that $\Pi_{\Phi_{\bar \l}} \neq \emptyset$, we just notice that the
following map
\begin{equation}\label{eq:ovPi}
  \ov{\pi}(\s,t) = t \Phi_{\bar \l} (\s); \qquad \quad \s \in \Sig_{k},
t \in [0,1] \quad ((\s,t) \in C_{k})
\end{equation}
belongs to $\Pi_{\Phi_{\bar \l}}$. Assuming by contradiction that $\a_\rho
\leq - \frac L2$, there would exist a map $\pi \in \Pi_{\Phi_{\bar \l}}$ with
$\sup_{\tilde{\s} \in C_{k}} II(\pi(\tilde{\s})) \leq - \frac 38 L$.
Then, since Proposition \ref{p:glob} applies with $\frac L4$,
writing $\tilde{\s} = (\s, t)$, with $\s \in \Sig_{k}$, the map
$$
  t \mapsto \Psi \circ \pi(\cdot,t)
$$
would be an homotopy in $\Sig_{k}$ between $\Psi \circ \Phi_{\bar \l}$ and a
constant map. But this is impossible since $\Sig_{k}$ is
non-contractible (see Lemma \ref{l:nonco}) and since $\Psi \circ
\Phi_{\bar \l}$ is homotopic to the identity, by Proposition \ref{eq:test}.
Therefore we deduce $\ov{\Pi}_{\Phi_{\bar \l}} > - \frac L2$.
\end{pf}

\begin{pfn} {\sc of Theorem \;$\ref{eq:theo}$}
We introduce a variant of the above minimax scheme,
following \cite{str} and \cite{djlw}. For $t$ close to $1$, we
consider the functional
\begin{eqnarray*}
  \text{II}_{t \rho_1, t \rho_2}(u) & = \frac{1}{2}\int_{\Sig}|\n_g u|^2dV_g-t\rho_{1}\log \int_{\Sig}e^{u-\bar u}dV_g-t\rho_2\log\int_{\Sig}e^{-u+\bar u}dV_g.
\end{eqnarray*}
Repeating the estimates of the previous sections, one easily checks
that the above minimax scheme applies uniformly for $t \in [ 1 -
t_0, 1 + t_0 ]$ with $t_0$ sufficiently small. More precisely, given
$L > 0$ as before, for $t_0$ sufficiently small we have
\ba\label{eq:minmaxrho} \nonumber
   \sup_{\pi \in \Pi_{\Phi_{\bar \l}}} \sup_{m \in \partial
   C_{k}} \text{II}_{t \rho_1,t \rho_2}(\pi(m))
   < - 2 L; \quad \a_{t \rho} := \inf_{\pi \in \Pi_{\Phi}}
  \; \sup_{m \in C_{k}} \text{II}_{t \rho_1, t \rho_2}(\pi(m)) > -
  \frac{L}{2}; \\ \hbox{ for every }  t \in [1 - t_0, 1 + t_0],
\ea where $\Pi_{\Phi_{\bar \l}}$ is defined in \eqref{eq:PiPi}.

Next we notice that for $t' \geq t$ there holds
$$
  \frac{\text{II}_{t \rho_1, t \rho_2}(u)}{t} - \frac{\text{II}_{t' \rho_1, t'
\rho_2}(u)}{t'} = \frac{1}{2}\left(\frac{1}{t}-\frac{1}{t^{'}}\right)\int_{\Sig}|\n_g u|^2dV_g \geq 0, \qquad u \in
H^1(\Sig) .
$$
Therefore it follows easily that also
$$
  \frac{\a_{t \rho}}{t} - \frac{\a_{t' \rho}}{t'} \geq 0,
$$
namely the function $t \mapsto \frac{\a_{t \rho}}{t}$ is
non-increasing, and hence is almost everywhere differentiable. Using
Struwe's monotonicity argument, see for example \cite{djlw}, one van
see that at the points where $\frac{\a_{t \rho}}{t}$ is
differentiable $\text{II}_{t \rho_1, t \rho_2}$ admits a bounded
Palais-Smale sequence at level $\a_{t \rho}$, which converges to a
critical point of $\text{II}_{t \rho_1, t \rho_2}$. Therefore, since the
points with differentiability fill densely the interval $[1 - t_0, 1
+ t_0]$, there exists $t_n \rightarrow 1$  and \;$u_n\in H^1(\Sig)$\;such that

\begin{equation}\label{eq:gtodak}
 -\D_gu_n=t_n\rho_1\left(\frac{e^{u_n}}{\int_{\Sig}e^{u_n}dV_g}-1\right)-t_n\rho_2\left(\frac{e^{-u_n}}{\int_{\Sig}e^{-u_n}dV_g}-1\right).
\end{equation}
At this stage , it is sufficient to  apply Proposition\;$\ref{eq:compactness}$\;to get a limit wich is a solution of \;$\eqref{eq:eqs}$. This conclude the proof.
\end{pfn}


\end{document}